\documentclass{amsart}
\usepackage{amssymb}
\usepackage{amsmath}

\newtheorem{theorem}{Theorem}[section]
\newtheorem{proposition}[theorem]{Proposition}
\newtheorem{lemma}[theorem]{Lemma}
\newtheorem{definition}[theorem]{Definition}

\numberwithin{equation}{section}
\begin{document}

\title[Degenerate $k$-Hessian Equations]{$C^{1, 1}$ Solution of the Dirichlet Problem for\\
Degenerate $k$-Hessian Equations}

\author{Qi Wang and Chao-Jiang Xu}
\address{School of Mathematics and Statistics, Wuhan University 430072,
Wuhan, P. R. China}
\email{qiwang88@whu.edu.cn  and  chjxu.math@whu.edu.cn}

\date{04/09/2013}



\maketitle

\begin{abstract}
In this paper, we prove the existence of $C^{1,1}$-solution to the Dirichlet problem for degenerate elliptic $k$-Hessian equations $S_{k}[u]=f$ under a condition which is weaker than the condition $f^{1/k}\in C^{1,1}(\bar\Omega)$.
\end{abstract}

\section{Introduction}
In this work, we study the following Dirichlet problem for the
$k$-Hessian
equation
\begin{equation}{\label{eq:o1}}
\left\{
\begin{array}{ll}
S_{k}[u]=f(x)&\indent\text{in}\,\,\,\,\Omega, \\
    u=\varphi(x)&\indent\text{on}\,\,\partial\Omega,
\end{array}
\right.
\end{equation}
where $\Omega$ is a bounded  domain in $\mathbb{R}^{n}$, $S_{k}[u]$ is
defined as follow
$$
S_{k}[u]=\sigma_{k}(\lambda),\indent k=1,\ldots,n,
$$
where $\lambda=(\lambda_{1},\lambda_{2},\ldots,\lambda_{n})$,
$\lambda_{i}$ is the
eigenvalue of the Hessian matrix $(D^{2}u)$, and
\begin{equation}\label{eq:ksigma}
\sigma_{k}(\lambda)=\sum_{1\leq i_{1}<\cdots<i_{k}\leq
n}\lambda_{i_{1}}\cdots\lambda_{i_{k}}
\end{equation}
is the $k$-th elementary polynomial. Note that the case $k=1$
corresponds to the Possion's
equation, while for $k=n$, that is the Monge-Amp\`{e}re equation $\det
D^{2}u=f$.

The nonlinear equation of \eqref{eq:o1} is referred to as non-degenerate when the
function $f$ is positive, it is degenerate elliptic if $f$
is non-negative and allowed to vanish somewhere in $\overline\Omega$.

The non-degenerate $k$-Hessian equations were firstly studied by
Caffarelli, Nirenberg and Spruck \cite{c2}. They proved that if
$f\in C^{1,1}(\overline{\Omega})$, $f>0$, $\partial\Omega$
and $\varphi$ were sufficiently smooth, (\ref{eq:o1}) had a unique
$C^{3,\,\alpha}$ $k$-admissible solution.
For the degenerate case, Ivochina, Trudinger and Wang \cite{I} studied
a class of fully nonlinear degenerate elliptic equations which depended
only on the eigenvalues of the Hessian matrix. This kind of equations
include the $k$-Hessian equations. They got the priori estimate with
the condition $f^{1/k}\in C^{1,1}(\overline\Omega)$. In particular,
their estimation of second order derivatives was
independent with $\inf_{\Omega} f$. Thus, the condition $f^{1/k}\in
C^{1,1}(\overline{\Omega})$ implied the existence of $C^{1,1}$-solutions
to the degenerate $k$-Hessian equations. Then, the regularity of the degenerate $k$-Hessian equations paused at $C^{1,1}$. For Monge-Amp\`{e}re equations, Hong, Huang and Wang \cite{H} gave a special condition to the smooth solution for the 2-dimensional Monge-Amp\`{e}re equation. We can find that even $f$ is analytic, the solution may be not in $C^{2}$ \cite{W1}.
 For degenerate $k$-Hessian equations, some papers concentrated on the convexity of the solutions \cite{M}.

In this work, we want to improve these results of $C^{1,1}$-regularity with a condition weaker
then $f^{1/k}\in C^{1,1}(\overline{\Omega})$. To state our results, we set
the following condition for the function $f$ which is the right hand side term of the
equations.

\smallskip
\noindent {
{\bf Condition (H) :}{\em \label{cond:f}
Assume that $f\in C^{1,1}(\overline{\Omega})$, $f\geq 0$ and  there exists a constant $C_0>0$ such that
$$
|Df(x)|\leq C_0f^{1-\frac 1k}(x)\quad\forall x\in{\Omega},
$$
and for any vector $\xi\in \mathbb{S}^{n-1}$,
$$
f(x)f_{\xi\xi}(x)-(1-\frac{1}{k})f^{2}_{\xi}(x)\geq
-C_0f^{2-\frac 1k}(x)\quad\forall x\in\overline{\Omega},
$$
where $f_{\xi}(x)=\frac{\partial f}{\partial\xi}(x),
f_{\xi\xi}(x)=\frac{\partial^2 f}{\partial\xi^2}(x)$.
}}

We will show that Condition {\bf(H)} is weaker than $f^{\frac 1k}\in C^{1,1}(\overline\Omega)$ in Section \ref{s2}. Indeed, for the
case of $3$-dimension we can give an example that $f\ge0$ is analytic and $f^{\frac 12}$ is only Lipshitz continuous, while $f$ satisfies Condition {\bf (H)}.

Our main result is stated as follow.

\begin{theorem}\label{th:exist}
Assume that $\Omega$ is a bounded $(k-1)$-convex domain in
$\mathbb{R}^{n}$ with $C^{3,1}$
boundary $\partial\Omega$, $f\geq 0$,  $f$ satisfies
Condition {\bf (H)}, and $\varphi\in C^{3,1}(\partial\Omega)$. Then the
Dirichlet problem
\eqref{eq:o1} has a unique $k$-admissible solution $u\in
C^{1,1}(\overline\Omega)$.
Moreover,
$$
\|u\|_{C^{1,1}(\overline{\Omega})}\leq C,
$$
where $C$ depends only on $n$, $k$, $\Omega$,
$\|f\|_{C^{1,1}(\overline{\Omega})}$,
$\|\varphi\|_{C^{3,1}(\partial\Omega)}$ and $C_0$. In
particular, $C$ is independent with
$\inf_{\Omega} f$.
\end{theorem}

We will recall the notions of $(k-1)$-convexity and $k$-admissibility
in Section 2.

In the paper \cite{G}, for the degenerate Monge-Amp\`{e}re equations, P. Guan introduced a
condition weaker then $f^{1/n}\in C^{1,1}(\overline{\Omega})$. So that our condition
{\bf (H)} is an extension of Guan's condition to $k$-Hessian equation.

In the following Section \ref{s2}, we will recall some definitions and some known results,
then we will give the sketch of the proof to the main theorem. Then, the rest of this paper (Section \ref{s3} to Section \ref{s5})
is to establish the uniform \`a priori estimates for the approximate solutions.

\bigskip

\section{Sketch of the proof to the Main Theorem}\label{s2}

In this section, we firstly recall  some definitions and known results. Then we will present the sketch of the proof of Theorem
\ref{th:exist}.

\subsection*{Preliminaries}

Firstly, we recall  some  definitions about the $k$-Hessian equations.

\begin{definition}[\cite{W}]\label{def:kadm}
We say a function $u\in C^{2}(\Omega)\cap C^{0}(\overline\Omega)$ is
k-admissible if
$$
\lambda(D^{2}u)\in \overline\Gamma_{k},
$$
where $\Gamma_{k}$ is an open symmetric convex cone in
$\mathbb{R}^{n}$, with vertex at the
origin, given by
$$
\Gamma_{k}=\{(\lambda_{1},\ldots,\lambda_{n})\in
\mathbb{R}^{n}:\,\,\sigma_{j}(\lambda)>0,\forall
j=1,\ldots,k\},
$$
where $\sigma_{j}(\lambda)$ is defined by \eqref{eq:ksigma}.
\end{definition}
The geometry condition for $\Omega\subset \mathbb{R}^{n}$ is (See \cite{c2}),
\begin{definition}\label{def:convex}
We say that $\Omega$ is $(k-1)$-convex if there exists a constant $c>0$, such that,
for any $x\in \partial\Omega$
$$
\sigma_{j}(\kappa)(x)\geq c>0\indent j=1,\ldots,k-1,
$$
where $\kappa=(\kappa_{1},\ldots,\kappa_{n-1})$, $\kappa_{i}(x)$ is the
principal curvature of $\partial\Omega$ at $x$. When $k=n$, it is the usual convexity.
\end{definition}
The weak solution to the $k$-Hessian equation is defined as follow.
\begin{definition}[\cite{T1}]\label{def:weakst}
A function $u\in C^{0}(\Omega)$ is called an admissible weak solution
of equation
 \eqref{eq:o1} in the domain $\Omega$, if there exists a sequence
$\{u_{m}\}\subset
C^{2}(\Omega)$ of $k-$admissible functions such that
$$
u_{m}\rightarrow u \indent \text{in}\,\, C^{0}(\Omega),\indent
S_{k}(u_{m})\rightarrow f
\indent \text{in}\,\,
L^{1}_{loc}(\Omega).
$$
\end{definition}

\smallskip
\subsection*{Condition {\bf (H)}}

For our Condition {\bf (H)}, by simple computation, we have
\begin{lemma}\label{lm:ch}
Assume that $f\in C^{1,1}(\overline\Omega), f\geq 0$ :
\begin{itemize}
 \item[(i)] if $f^{1/k}\in C^{1,1}(\overline\Omega)$, then, $f$ satisfies {\textup{Condition}} {\bf (H)}.
 \item[(ii)] if $f\geq \delta_{0}>0$ and $f$ satisfies {\textup{Condition}} {\bf(H)},
 then, $f^{1/k}\in C^{1,1}(\overline\Omega)$;
 \item[(iii)] if $f$ satisfies {\textup{Condition}} {\bf(H)}, then for any small $\epsilon>0$, $f+\epsilon$ satisfies
 {\textup{Condition}} {\bf(H)} with the same constant $C_0$.
\end{itemize}
\end{lemma}

The following interesting example shows the signification of our Condition {\bf(H)}.

\noindent {\bf Example.}
Let $\Omega=B_{1}(0)=\{(x, y, z)\in \mathbb{R}^{3};\,\, |(x, y, z)|<1\}$ and
$$u=(x^{2}+y^{2}+z^{2})^{3/2}-1.$$
Then, we have
\begin{equation*}
\left\{
\begin{array}{ll}
S_{2}[u]=45(x^{2}+y^{2}+z^{2})\qquad&\,\text{in}\,\,\,\Omega, \\
    u=0&\text{on}\,\partial\Omega.
\end{array}\right.
\end{equation*}
Remind that $u\in C^{2,1}(\overline{\Omega})$ and $u\notin
C^{3}(\overline{\Omega})$. In particular, $f=45(x^{2}+y^{2}+z^{2})$ is analytic and $f^{1/2}$ is only
Lipschtz continuous near the origin. However, $f$ satisfies Condition {\bf (H)}.
Indeed, $f$ is radially symmetric. Then, we may choose $\xi=(1,0,\ldots,0)$,
for any positive constant $C_0$,
\begin{align*}
f(x)f_{\xi\xi}(x)-(1-\frac{1}{2})f^{2}_{\xi}(x)&=
\frac{90^{2}}{2}(x^{2}+y^{2}+z^{2})-\frac{90^{2}}{2}x^{2}
\geq 0\\
&\geq -C_0 \big(45(x^{2}+y^{2}+z^{2})\big)^{3/2}=-C_0f^{2-1/2}.
\end{align*}

\smallskip
\subsection*{Sketch of the proof to Theorem \ref{th:exist}}
The proof of the Theorem \ref{th:exist} is standard by using a non-degenerate approximation and
a uniform \`a priori estimate of approximate solution.

For the non-degenerate equation, we have the following existence result of Caffarelli, Nirenberg, Spruck\cite{c2} and Trudinger \cite{T}.
\begin{proposition}\label{th:c2}
Assume that $\Omega$ is a bounded $(k-1)$-convex domain in
$\mathbb{R}^{n}$ with $C^{3,1}$ boundary $\partial\Omega$,
 $f^{1/k}\in
C^{1,1}(\overline{\Omega})$ and $f\geq \delta_{0}$ for some positive constant $\delta_{0}$,
$\varphi\in C^{3,1}(\partial\Omega)$. Then
the Dirichlet problem \eqref{eq:o1} has a unique $k$-admissible solution
$u\in C^{3,\alpha}(\overline\Omega)$. Moreover,
$$\|u\|_{C^{3,\alpha}}\leq C,$$
where $\alpha\in (0,1)$, $C$ depends only on $n$, $k$, $\alpha$,
$\delta_{0}$, $\Omega$, $\|\varphi\|_{C^{3,1}(\overline{\Omega})}$, and
$\|f\|_{C^{1,1}(\overline{\Omega})}$.
\end{proposition}

Now we study the following approximate problem with $\vartheta>0$,
\begin{equation}\label{eq:1}
\left\{
\begin{array}{ll}
S_{k}[u^{\vartheta}]=f+\vartheta &\indent\indent \text{in}\,\,\,\,\Omega,\\
u^{\vartheta}=\varphi &\indent\indent\text{on}\,\partial\Omega.
\end{array}\right.
\end{equation}
We will prove the following uniform \`a priori estimate.

\begin{theorem}[\`a priori estimate]\label{th}
Assume that $\Omega$ is a bounded $(k-1)$-convex domain in
$\mathbb{R}^{n}$ with $C^{3,1}$
boundary $\partial\Omega$,
$f$ satisfies
Condition {\bf (H)}, $\varphi\in C^{3,1}(\partial\Omega)$. Let $u^{\vartheta}\in
C^{3, 1}(\overline{\Omega})$ be a $k$-admissible solution to the
Dirichlet problem
\eqref{eq:1}. Then, we have the following \`a priori estimate
$$
\|u^{\vartheta}\|_{C^{1,1}(\overline{\Omega})}\leq C,
$$
where $C$ depends only on $n$, $k$, $\Omega$,
$\|f\|_{C^{1,1}(\overline{\Omega})}$,
$\|\varphi\|_{C^{2,1}(\partial\Omega)}$ and $C_0$. In
particular, $C$ is independent with
$\vartheta$.
\end{theorem}

The proof of the theorem above will be our main task for the rest of
this paper. Now we explain how to prove Theorem \ref{th:exist} through
Theorem \ref{th}.

\smallskip
\noindent
\begin{proof}[{\bf Proof of Theorem \ref{th:exist}}]
As we have said
previously, we complete the proof of the existence theorem by an
approximation of non-degenerate problems, (See \cite{G} and \cite{GTW}).
Since $f$ satisfies Condition {\bf(H)}, $f+\vartheta\geq\vartheta>0$. By Lemma \ref{lm:ch}, $(f+\vartheta)^{1/k}\in C^{1,1}(\overline\Omega)$. Thus, by Proposition \ref{th:c2}, Theorem \ref{th} and the continuity method in \cite{E} and \cite{GT}, the equation \eqref{eq:1} has a $k$-admissible solution $u^\vartheta$ which belongs to $C^{3,\alpha}(\overline\Omega)$ and
satisfies the following  estimate
$$
\|u^{\vartheta}\|_{C^{1,1}}\leq C
$$
with $C$ depends only on $n$, $k$, $\Omega$, $\|f\|_{C^{1,1}(\overline{\Omega})}$,
$\|\varphi\|_{C^{2,1}(\partial\Omega)}$ and $C_0$, in particular, $C$ is independent with $\vartheta$.
Then, the Azel\`{a}-Ascoli Theorem implies that the Dirichlet problem (\ref{eq:o1}) admits a $k$-admissible weak
solution $u\in C^{1,1}(\overline\Omega)$. The uniqueness from the comparison principles of fully
nonlinear degenerate elliptic equations.
\end{proof}

\bigskip
\section{Interior \`a Priori estimates}\label{s3}

Without any difficulties, we can get the $L^{\infty}$ estimate
and the gradient estimate of the approximate solution similar to the papers \cite{c2} and
\cite{W}.

\begin{proposition}\label{lm:222}
Assume that $\Omega$ is a bounded $(k-1)$-convex domain in
$\mathbb{R}^{n}$ with $C^{3,1}$ boundary. $f$ satisfies Condition {\bf(H)},
$\varphi\in C^{3,1}(\partial\Omega)$. Let $u^{\vartheta}$ be a  $k$-admissible
solution to the Dirichlet problem \eqref{eq:1}. Then we have
\begin{equation*}
\|u^{\vartheta}\|_{C^{1}(\overline\Omega)}\leq C,
\end{equation*}
where $C$ depends only on $\|\varphi\|_{C^{1}(\partial\Omega)}$,
$\|f\|_{C^{1}(\overline{\Omega})}$, $C_0$ and $\Omega$, in particular, $C$ is independent with $\vartheta$.
\end{proposition}

Now we consider the second order derivative estimates of $u^{\vartheta}$ which suffice to
prove the following one side estimate,
\begin{equation}\label{secesti}
u^{\vartheta}_{\xi\xi}(x)\leq C\qquad \forall (x,\xi)\in\bar\Omega\times\mathbb{S}^{n-1}.
\end{equation}
Since we have already known that $\Delta u^{\vartheta}\geq 0$ by the definition of the $k$-admissibility,
then by a rotation, we obtain,
$$
|u^{\vartheta}_{\xi\xi}(x)|\leq nC.
$$
Thus, if we proved (\ref{secesti}), we could finish the whole proof of
Theorem \ref{th}.

Now we give  a relationship between the boundary estimate and the interior
estimate. That is
\begin{lemma}\label{th1}
Assume that $\Omega$ is a $(k-1)$-convex domain in $\mathbb{R}^{n}$
with $C^{3,1}$
boundary. $f$ satisfies
the condition
{\bf (H)}, $\varphi\in C^{3,1}(\partial\Omega)$. Let $u^{\vartheta}$ be a
$k$-admissible solution
to the Dirichlet problem \eqref{eq:1}. Then we have
\begin{equation}\label{eq:rela}
\sup_{\Omega}u^{\vartheta}_{\eta\eta}\leq C+\sup_{\partial\Omega} u^{\vartheta}_{\eta\eta}
\end{equation}
for any unit vector $\eta\in \mathbb{S}^{n-1}$, where $C$ depends only on
$\Omega$ and $C_0$, in particular, $C$ is independent with $\vartheta$.
\end{lemma}

\begin{proof}
Denote
\begin{equation}\label{eq:denoF}
F[D^2 u^{\vartheta}]=\big(S_k(u^\vartheta)\big)^{1/k}=(f+\vartheta)^{1/k}.
\end{equation}
Then, differentiating (\ref{eq:denoF}) in direction
$\eta\in\mathbb{S}^{n-1}$ twice, one can
verify that
$$
F^{ij}u^{\vartheta}_{ij\eta\eta}+F^{ij,st}u^{\vartheta}_{ij\eta}u^{\vartheta}_{st\eta}=\frac{1-k}{k^{2}}
(f+\vartheta)^{1/k-2}f_{\eta}^{2}+\frac{1}{k}
(f+\vartheta)^{1/k-1}f_{\eta\eta},
$$
where
$$
F^{ij}=F^{u^{\vartheta}_{ij}}=\frac{\partial F[u^{\vartheta}]}{\partial
u^{\vartheta}_{ij}},\qquad F^{ij,st}=F^{u^{\vartheta}_{ij},u^{\vartheta}_{st}}=\frac{\partial^{2}F[u^{\vartheta}]}{\partial
u^{\vartheta}_{ij}\partial u^{\vartheta}_{st}}\,.
$$

From \cite{c2} and \cite{W}, we have known that $F$ is a concave operator,
thus,
$$F^{ij}u^{\vartheta}_{ij\eta\eta}\geq
\frac{1-k}{k^{2}}(f+\vartheta)^{1/k-2}f^{2}_{\eta}+\frac{1}{k}(f+\vartheta)^{1/k-1}f_{\eta\eta}.$$
That is
\begin{equation*}
\mathfrak{L}u^{\vartheta}_{\eta\eta}\geq \frac{1-k}{k}f^{2}_{\eta}+(f+\vartheta)f_{\eta\eta},
\end{equation*}
where $\mathfrak{L}=(f+\vartheta)\cdot S^{ij}_{k}[u^{\vartheta}]\partial_{i}\partial_{j}$.
Next, let $w=\frac{1}{2}a|x|^{2}$, we have, by Maclaurin's inequality
$$\mathfrak{L}w=a(n+1-k)\sigma_{k-1}[\lambda(D^{2}u^{\vartheta})]\cdot (f+\vartheta)\geq a
\cdot C_{n,k}(f+\vartheta)^{2-1/k},$$
where $C_{n,k}$ is a constant depends only on $n$ and $k$.
By Condition {\bf (H)}, we choose $a$ so large that
$$a\cdot C_{n,k}\geq C_{0}, $$
then,
$$\mathfrak{L}[w+u^{\vartheta}_{\eta\eta}]\geq
C_{0}(f+\vartheta)^{2-1/k}+\frac{1-k}{k}f^{2}_{\eta}+(f+\vartheta)f_{\eta\eta} \geq 0$$
holds.  By the classical weak maximum principle,
$$\sup_{\Omega}(w+u^{\vartheta}_{\eta\eta})\leq
\sup_{\partial\Omega}(w+u^{\vartheta}_{\eta\eta}).$$
That is (\ref{eq:rela}), proof is done.
\end{proof}

\subsection*{Tangential and mixed second order derivative estimates}
By Lemma \ref{th1}, we reduce the estimate to the boundary.
For any given boundary point $x_{0}\in \partial\Omega$, by a
translation and rotation of
coordinates, we assume that $x_{0}$ is the origin and locally
$\partial\Omega$ is given by
$$x_{n}=\rho(x'),x'=(x_{1},\ldots,x_{n-1}),$$
such that $D\rho(x_0)=0$. Exactly as which in \cite{c2} we have the following result.

\begin{proposition}\label{lm:3}
Assume that $\Omega$, $f$ and $\varphi\in C^{3,1}(\partial\Omega)$
satisfy the conditions of which in Theorem \ref{th}. Let $u^{\vartheta}$ be a $k$-admissible
solution to the Dirichlet problem
\eqref{eq:1}, $x_0\in\partial\Omega$. Then we have
$$
|u^{\vartheta}_{ij}(x_0)|\leq C,\indent \, |u^{\vartheta}_{in}(x_0)|\leq C,\indent i, j=1,
\ldots, n-1,
$$
where the constant $C$ depends only on $\Omega$, $C_{0}$,
$\|\varphi\|_{C^{2,1}(\partial\Omega)}$ and
$\|f\|_{C^{1,1}(\overline{\Omega})}$, in particular, $C$ is independent
with $\vartheta$.
\end{proposition}

\bigskip

\section{The Weakly Interior Estimate}\label{s4}
Since we have completed the double tangent derivative estimate and the mix type derivative estimate by
Proposition \ref{lm:3}, we just need to study the double normal derivative estimate. We can assume that
\begin{equation}\label{amax}
\max\{\sup_{\partial\Omega}|u^{\vartheta}|,\sup_{\partial\Omega}|Du^{\vartheta}|,\sup_{\partial\Omega}
D^{2}u^{\vartheta}\}=\sup_{x\in\partial \Omega} u^{\vartheta}_{\gamma\gamma}(x),
\end{equation}
where $\gamma\in\mathbb{S}^{n-1}$ is normal vector of $\partial\Omega$ at $x\in\partial\Omega$. If \eqref{amax} did not hold, we should have finished the proof by Proposition \ref{lm:3}. Besides,
we also assume $\sup_{x\in\partial\Omega} u^\vartheta_{\gamma\gamma}(x)>0$.

The double normal derivative estimate will be established by two steps.
The first step is the following weakly interior estimate.

\begin{lemma}\label{lm:weakin}
Assume that function $u^{\vartheta}\in C^{3,1}(\overline\Omega)$
satisfies
\eqref{eq:1}, $\Omega$ is $(k-1)$-convex, $\partial\Omega\in C^{3,1}$,
$\varphi\in
C^{3,1}(\partial\Omega)$,
$f$ satisfies
 Condition {\bf (H)}. Then we have a weakly interior estimate. That is,
for some
sufficiently small constant $\varepsilon$ and some $\delta>0$, there
exists a constant
$C_{\varepsilon,\delta}$ such that, for any $\xi\in \mathbb{R}^{n}$,
\begin{equation}\label{eq:weakin}
\sup_{x\in\Omega_{\delta}} u^{\vartheta}_{\xi\xi}(x)\leq
\varepsilon\sup_{x\in\partial\Omega}u^{\vartheta}_{\gamma\gamma}(x)+C_{\varepsilon,\delta}|\xi|^{2},
\end{equation}
where $\Omega_{\delta}=\{x\in \Omega|\textup{dist}(x,\partial\Omega)\geq\delta\}$,
$\gamma$ is the unit
inner normal vector, $C_{\varepsilon,\delta}$ depends only on $n$,
$\varphi$, $\partial\Omega$, $C_{0}$ and
$\|f\|_{C^{1,1}(\overline\Omega)}$, in particular,
$C_{\varepsilon,\delta}$ is independent with $\vartheta$.
\end{lemma}

To prove the proposition above, we will use the following maximum principle.
\begin{proposition}[{\cite{I}}]\label{lm:nmp}
Let $A=A(x,\tilde\xi)$ be a $(n+N)\times (n+N)$ matrix which is positive
definite,
${\bf{b}}={\bf{b}}(x,\tilde\xi)$ be a $(n+N)$-dimensional vector in $Q$, where
$Q=:\Omega\times
\mathbb{R}^{N}$, $x\in \Omega$, $\tilde\xi\in\mathbb{R}^{N}$, such that
\begin{equation}\label{defofL}
\mathcal{L}[w]=:tr(AD^{2}w)+<\mathbf{b},Dw>
\end{equation}
is an elliptic operator in $Q$, where $tr(AD^{2}w)$ is the trace of
matrix $(AD^{2}w)$,
$<,>$ is the inner product. Assume that $g,h\in C^{2}(\overline Q)$,
$h>0$, and $g,h$ are
$p$-homogeneous in $\tilde\xi$ for some $p>1$. If there exists positive
constants $\mu,\nu$ such
that
\begin{equation}\label{eq:nmp1}
\mathcal{L}[g]\geq-\mu|\tilde\xi|^{p}\indent \forall (x,\tilde\xi)\in \Omega\times \mathbb{R}^{N},
\end{equation}
\begin{equation}\label{eq:nmp2}
\mathcal{L}[h]\leq -\nu|\tilde\xi|^{p}\indent \forall (x,\tilde\xi)\in \Omega\times \mathbb{R}^{N}.
\end{equation}
Then, we have
$$
\sup_{\Omega\times\{|\tilde\xi|=1\}}\frac{g}{h}\leq \frac{\mu}{\nu}+\sup_{\partial\Omega\times\{|\tilde\xi|=1\}}\frac{g}{h}.
$$
\end{proposition}

Next, we construct the operator
$\mathcal{L}$ which has the form \eqref{defofL}, the functions $g$ and $h$ which are applicable to  Proposition \ref{lm:nmp}. Denote $\tilde{\xi}=(\xi,\xi_{0})$, $\xi=(\xi_{1},\xi_{2},\ldots,\xi_{n})\in \mathbb{R}^{n}$, $\xi_{0}\in \mathbb{R}$. That is $Q\ni(x,\tilde{\xi})=\Omega\times\mathbb{R}^{N}=: \Omega\times\mathbb{R}^{n+1}$.  The matrix $A$
and vector $\mathbf b$ are given by the following formula:
\begin{gather*}
A=\begin{pmatrix}
\mathcal{A}&r\mathcal{A}&\mathcal{A}q\\
r\mathcal{A}&r^{2}\mathcal{A}&r\mathcal{A}q\\
q'\mathcal{A}&rq'\mathcal{A}&q'\mathcal{A}q
\end{pmatrix}
,
\mathbf{b}=\begin{pmatrix}
0\\
-2\mathcal{A}q\\
0
\end{pmatrix}
,
\end{gather*}
where
$$
\mathcal{A}=\frac{1}{tr(\mathcal{F})}\mathcal{F}=\big( a^{j k}\big)\,,(F^{ij})_{n\times n}=\mathcal{F},
$$
and
\begin{equation*}
\begin{split}
r&=\frac{\alpha\psi_{\xi}}{\psi+\beta},\\
q&=\frac{1}{\psi+\beta}\left(\frac{\xi}{4}+\frac{\alpha\psi_{\xi}}{\psi+\beta}D\psi\right),
\end{split}
\end{equation*}
$0<\alpha<1/2,0<\beta\ll1$ and $0<\hat{\varepsilon}\ll\beta^{4}$. $\psi$ is a special auxiliary function (see Section 2 of \cite{c2} and the Lemma 3.3 of \cite{I0} )
 which is non-negative and $C^{\infty}$ on $\overline\Omega$ such that
\begin{equation}\label{eq:psipro}
|D\psi|\geq n\indent \text{on}\,\partial\Omega,\indent
tr(\mathcal{A}D^{2}\psi)\leq
-n\indent \text{in}\,\Omega\,.
\end{equation}

The functions
$g$ and $h$ are defined as follows
\begin{equation}\label{eq:defg}
g=u^{\vartheta}_{\xi\xi}+2\xi_{0}u^{\vartheta}_{\xi}+\xi^{2}_{0}u^{\vartheta},
\end{equation}
\begin{equation}\label{eq:defh}
h=\frac{1}{4\alpha}(\psi+\beta)^{1-\alpha}|\xi|^{2}+\frac{\psi^{2}_{\xi}}{(\psi+\beta)^{\alpha}}+
\hat{\varepsilon}(\psi+\beta)
\xi_{0}^{2}.
\end{equation}
Immediately,
$g$ and $h$ are $2$-homogeneous in $\tilde \xi$.
In order to apply Proposition \ref{lm:nmp}, we need to confirm
assumptions (\ref{eq:nmp1}) and (\ref{eq:nmp2}). Those are the following two lemmas. For \eqref{eq:nmp1}, we have

\begin{lemma}\label{lm:nmp1}
To the function $g$ defined as \eqref{eq:defg}, we have \eqref{eq:nmp1}
holds, that is
$$\mathcal{L}[g]\geq-\mu|\tilde\xi|^{2}\indent \forall (x,\tilde{\xi})\in \Omega\times\mathbb{R}^{N},$$
where $\mu$ can be $\beta^{-5}$.
\end{lemma}
\begin{proof}
We have
\begin{equation*}
\begin{split}
\mathcal{L}[g]&=a^{ij}[u^{\vartheta}_{ij\xi\xi}+(4r+2\xi_{0})u^{\vartheta}_{\xi
ij}+(2r^{2}+4r\xi_{0}+\xi_{0}^{2})u^{\vartheta}_{ij}+4ru^{\vartheta}_{i}q_{j}
+2q_{i}q_{j}u^{\vartheta}]\\
&=a^{ij}(\tilde{I}_{1}+\tilde{I}_{2}+\tilde{I}_{3}+\tilde{I}_{4}+\tilde{I}_{5}).
\end{split}
\end{equation*}
For each term of $\mathcal{L}[g]$, we have

$$a^{ij}\tilde{I}_{1}=tr(\mathcal{A}D^{2}u^{\vartheta}_{\xi\xi})\geq
-K|\xi|^{2},$$
\begin{equation*}
|a^{ij}\tilde{I}_{2}|\leq|4rK+2K\xi_{0}||\xi|\leq K(|\xi|^{2}+\xi_{0}^{2})+\frac{C(\alpha,K,c_{0})}{\beta}|\xi|^{2},
\end{equation*}
where $$c_{0}=\|\psi\|_{C^{3}(\Omega)},\indent K=\max\{\|f\|^{1/k}_{C^{1,1}(\overline\Omega)},C_{0}\}.$$
For $a^{ij}\tilde{I}_{3}$, since what we need is the lower bound, we only need to estimate
$4r\xi_{0}u^{\vartheta}_{ij}a^{ij}$, that is
$$|4r\xi_{0}u^{\vartheta}_{ij}a^{ij}|\leq\frac{C(\alpha,K,c_{0})}{\beta}(|\xi|^{2}+\xi_{0}^{2}).$$
For the left terms, one can verify that
\begin{equation*}
\begin{split}
|a^{ij}\tilde{I}_{4}|&=\left|4a^{ij}u^{\vartheta}_{i}\frac{\alpha\psi_{\xi}}{(\psi+\beta)^{2}}\left(\frac{\xi_{j}}{4}+
\frac{\alpha\psi_{j}\psi_{\xi}}{\psi+\beta}\right)\right|\\
&\leq\frac{C(n,\alpha,c_{0},\|u\|_{L^{\infty}})}{\beta^{3}}|\xi|^{2},\\
\end{split}
\end{equation*}
and
\begin{equation*}
\begin{split}
|a^{ij}\tilde{I}_{5}|&\leq|2u^{\vartheta}\cdot q^{2}_{i}|\leq
\frac{2\|u^{\vartheta}\|_{L^{\infty}}}{\beta^{2}}\left(\frac{1}{16}+\frac{n^{2}\alpha^{2}}{\beta^{2}}c_{0}^{4}+
\frac{n\alpha c_{0}^{2}}{2\beta}\right)|\xi|^{2}\\
&\leq
\frac{C(n,\alpha,c_{0},\|u\|_{L^{\infty}})}{\beta^{4}}|\xi|^{2},
\end{split}
\end{equation*}
Thus,
\begin{equation*}
\begin{split}
\mathcal{L}[g]&\geq-(2K+\frac{C(\alpha,K,c_{0})}{\beta}+\frac{C(n,\alpha,c_{0},\|u\|_{L^{\infty}})}{\beta^{3}}
+\frac{C(n,\alpha,c_{0},\|u\|_{L^{\infty}})}{\beta^{4}})|\xi|^{2}\\
&-(K+\frac{C(\alpha,K,c_{0})}{\beta})\xi_{0}^{2}.
\end{split}
\end{equation*}
Since $\beta$ is sufficiently small, we may choose
\begin{equation}\label{eq:defmu}
\mu=\beta^{-5}.
\end{equation}

Thus, we have
$$
\mathcal{L}[g]\geq -\mu(|\xi|^{2}+\xi_{0}^{2})=-\mu|\tilde{\xi}|^{2}.
$$
Proof is done.
\end{proof}

Then, we state \eqref{eq:nmp2} as follow.
\begin{lemma}\label{lm:nmp2}
To the function $h$ defined as \eqref{eq:defh}, we have
$$\mathcal{L}[h]\leq-\nu|\tilde{\xi}|^{2}\indent \forall (x,\tilde\xi)\in \Omega\times \mathbb{R}^{N},$$
where
$$\nu=\min\left\{\frac{1}{16\alpha(K_{0}\delta+\beta)^{\alpha}},\frac{\hat{\varepsilon}}{2}\right\}.$$
\end{lemma}
\begin{proof}
We divide $h$ into two parts: $\hat{h}$ and $h-\hat{h}$. $\hat{h}$ is
defined as
$$
\hat{h}=\frac{1}{4\alpha}(\psi+\beta)^{1-\alpha}\xi^{2}+\frac{\psi^{2}_{\xi}}{(\psi+\beta)^{\alpha}}.
$$

We have
\begin{equation*}
\begin{split}
\mathcal{L}[\hat{h}]&=tr(A\cdot D^{2}\hat{h})+<\mathbf{b},D\hat{h}>\\
&=\frac{1}{N_{0}(\psi+\beta)^{\alpha+1}}\left\{F^{ij}(\psi+\beta)I_{1}+F^{ij}I_{2}
+\frac{F^{ij}}{(\psi+\beta)}I_{3}\right.\\
&\left.+F^{ij}(\psi+\beta)(-\frac{\xi_{i}\xi_{j}}{4\alpha})+F^{ij}
(-\alpha-1)\psi_{i}\psi_{\xi}\xi_{j}\right\},
\end{split}
\end{equation*}
where
\begin{equation*}
\begin{split}
I_{1}&=\frac{1-\alpha}{4\alpha}\psi_{ij}|\xi|^{2}+2\psi_{\xi
i}\cdot\psi_{\xi
j}+2\psi_{\xi}\psi_{\xi ij},\\
I_{2}&=3\alpha\psi_{ij}\cdot\psi_{\xi}^{2}-\frac{1-\alpha}{4}\psi_{i}\psi_{j}
|\xi|^{2}+\frac{\alpha}{2}\psi_{\xi}^{2}\delta_{ij},\\
I_{3}&=-3\alpha\psi_{i}\psi_{j}\psi_{\xi}^{2}-\alpha^{2}\psi_{i}\psi_{j}\psi_{\xi}^{2}.
\end{split}
\end{equation*}
To $F^{ij}I_{1}$, we have, by (\ref{eq:psipro}),
$$\frac{1-\alpha}{4\alpha}F^{ij}\psi_{ij}|\xi|^{2}\leq-\frac{1-\alpha}{4\alpha}N_{0}|\xi|^{2},\indent 2F^{ij}\psi_{\xi i}\psi_{\xi j}\leq 2N_{0}c^{2}_{0}|\xi|^{2},$$
$$2F^{ij}\psi_{\xi}\psi_{\xi ij}\leq
2n^{2}N_{0}\|D\psi\|_{L^{\infty}}\|D^{3}\psi\|_{L^{\infty}}|\xi|^{2}\leq
 2n^{2}N_{0}c_{0}^{2}|\xi|^{2}.$$
Thus, we have
\begin{equation}\label{eq:I1}
F^{ij}I_{1}\leq[-\frac{1-\alpha}{4\alpha}+(2n^{2}+2)c_{0}^{2}]N_{0}|\xi|^{2}.
\end{equation}
To $F^{ij}I_{2}$, one can verify that, by \eqref{eq:psipro}
$$3\alpha
F^{ij}\psi_{ij}\psi_{\xi}^{2}+\frac{\alpha}{2}\psi_{\xi}^{2}\delta_{ij}F^{ij}\leq
(-3n\alpha\psi_{\xi}^{2}+\frac{n\alpha}{2}\psi_{\xi}^{2})N_{0}\leq
0,$$
$$-\frac{1-\alpha}{4}F^{ij}\psi_{i}\psi_{j}|\xi|^{2}\leq0,$$
Thus, we have
\begin{equation}\label{eq:I2}
F^{ij}I_{2}\leq 0.
\end{equation}
For the left terms, we have
\begin{equation}\label{eq:left}
-(3\alpha+\alpha^{2})\frac{F^{ij}}{\psi+\beta}\psi_{i}\psi_{j}\psi_{\xi}^{2}-F^{ij}(\psi+\beta)\frac{\xi_{i}\xi_{j}}{4\alpha}-(1+\alpha)F^{ij}\psi_{i}\psi_{xi}
\xi_{j}\leq 0.
\end{equation}
In sum of the computation (\ref{eq:I1}), (\ref{eq:I2}) and \eqref{eq:left},
\begin{equation*}
\mathcal{L}[\hat{h}]\leq
\frac{1}{(\psi+\beta)^{\alpha}}\left[-\frac{1-\alpha}{4\alpha}+(2n^{2}+2)c_{0}^{2}\right]|\xi|^{2}.
\end{equation*}
We can find a constant $\alpha_{0}$ which depends on $c_{0}$ such that,
for
$0<\alpha\leq\alpha_{0}$
$$-\frac{1-\alpha}{4\alpha}+(n^{2}+2)c_{0}^{2}\leq
-\frac{1}{8\alpha}.$$
Thus, we have
\begin{equation}\label{eq:lhesti}
\mathcal{L}[\hat{h}]\leq-\frac{|\xi|^{2}}{8\alpha(K_{0}\delta+\beta)^{\alpha}}.
\end{equation}
For the left terms, we have
\begin{equation*}
\mathcal{L}[\hat{\varepsilon}(\psi+\beta)\xi_{0}^{2}]
\leq\hat{\varepsilon}
[-\xi^{2}_{0}+\frac{C(n,\alpha,c_{0})}{\beta^{2}m}\xi^{2}_{0}+\frac{C(n,\alpha,c_{0})
m}{\beta^{2}}|\xi|^{2}],
\end{equation*}
where $m$ is sufficiently large.

We choose $m$ so large that
$$-1+\frac{C(n,\alpha,c_{0})}{\beta^{2}m}\leq-\frac{1}{2}.$$
Still, we can choose $\hat{\varepsilon}$ sufficiently small such that
$$\hat{\varepsilon}\frac{C(n,\alpha,c_{0})m}{\beta^{2}}\leq
\frac{1}{16\alpha(K_{0}\delta+\beta)^{\alpha}}.$$
Thus, we claim, by (\ref{eq:lhesti}),
$$\mathcal{L}[h]\leq
-\frac{1}{16\alpha(K_{0}\delta+\beta)^{\alpha}}|\xi|^{2}-\frac{\hat{\varepsilon}}{2}\xi^{2}_{0}.$$
We denote
\begin{equation}\label{eq:defnu}
\nu=\min\left\{\frac{1}{16\alpha(K_{0}\delta+\beta)^{\alpha}},\frac{\hat{\varepsilon}}{2}\right\},
\end{equation}
we may choose $\nu=\beta^{3}$, hence
\begin{equation*}
\mathcal{L}[h]\leq-\nu(|\xi|^{2}+\xi_{0}^{2}).
\end{equation*}
\end{proof}

By Lemma \ref{lm:nmp1} and Lemma \ref{lm:nmp2}, we can establish the
proof of our weakly interior estimate.

\begin{proof}[{\bf{Proof of Lemma \ref{lm:weakin}}.}]
By the virtue of Proposition \ref{lm:nmp}
\begin{equation}\label{eq:nmp3'}
\sup_{\Omega\times\{|\tilde{\xi}|=1\}}\frac{g}{h}\leq \sup_{\partial\Omega\times\{|\tilde\xi|=1\}}\frac{g}{h}+\frac{\mu}{\nu}.
\end{equation}
As we have proved, $\frac{\mu}{\nu}$ depends on $\|f\|_{C^{1,1}}$,
$C_{0}$,
$\|\varphi\|_{L^{\infty}},\|\psi\|_{C^{3}},\delta$ and
$\|\rho\|_{C^{2,1}}$.

Next, we estimate $\sup_{\partial\Omega\times\{|\tilde\xi|=1\}}(g/h)$.

Let $z=(x,\xi,\xi_{0})\in\partial \Omega\times\{|\tilde\xi|=1\},\,|\xi|^{2}+\xi_{0}^{2}=1$ such
that
$$\frac{g}{h}(z)\geq\frac{1}{2}\sup_{\partial\Omega\times\{|\tilde\xi|=1\}}\frac{g}{h}.$$
Denote $\tilde\theta=<\gamma,\xi>=|\xi|\cos\theta$, where $\gamma$ is the unit normal of
$\partial\Omega$ at $x$.
Thus, at $z$ we have, there exists a unit tangent vector $\tau$ such that
$$\xi=\gamma\cdot|\xi|\cos\theta+\tau\cdot|\xi|\sin\theta.$$
Then,
\begin{equation*}
\begin{split}
u^{\vartheta}_{\xi\xi}&=u^{\vartheta}_{\gamma\gamma}|\xi|^{2}\cos^{2}\theta+u^{\vartheta}_{\tau\tau}|\xi|^{2}\sin^{2}\theta+
2u^{\vartheta}_{\tau\gamma}|\xi|^{2}
\sin\theta\cdot\cos\theta\\
&\leq2\sup_{\partial\Omega}u^{\vartheta}_{\gamma\gamma}\cdot\tilde\theta^{2}+2\|u^{\vartheta}\|_{C^{1}}(1-\tilde\theta^{2}),
\end{split}
\end{equation*}

Thus, we have
$$g\leq
2\tilde\theta^{2}\sup_{\partial\Omega}u^{\vartheta}_{\gamma\gamma}+n\|u^{\vartheta}\|_{C^{1}}(1-\theta^{2}).$$
Still,  at $z$, $\psi=0$ and $D\psi$ is in parallel with $\gamma$,  if
$\hat{\varepsilon}\leq\beta^{1-\alpha}/4\alpha$,
\begin{equation*}
\begin{split}
h&\geq \hat{\varepsilon}|\xi|^{2}+\theta^{2}\cdot
\beta^{-\alpha}+\varepsilon\beta\xi^{2}_{0}\\
&\geq \hat{\varepsilon}(1-\theta^{2})+\theta^{2}\cdot \beta^{-\alpha}.
\end{split}
\end{equation*}
$\hat{\varepsilon}$ may change from line to line, but all of them are
smaller than
$\beta^{1-\alpha}/4\alpha$, we  can choose
$\hat{\varepsilon}=\beta^{3}$.
Hence we can obtain
$$\sup_{\partial\Omega\times\{|\tilde\xi|=1\}}\frac{g}{h}\leq
\beta^{\alpha}\sup_{\partial\Omega}u^{\vartheta}_{\gamma\gamma}+\frac{n\|u^{\vartheta}\|_{L^{\infty}}}{\hat{\varepsilon}}\leq
\beta^{\alpha}\cdot
\sup_{\partial\Omega}u^{\vartheta}_{\gamma\gamma}+\hat{c}_{0}\beta^{-3},$$
where $\hat{c}_{0}$ depends on $\|u^{\vartheta}\|_{L^{\infty}}$.
Then, by the virtue of Proposition \ref{lm:nmp} and (\ref{eq:nmp3'}), for any
$z\in \Omega\times\{|\tilde\xi|=1\}$
\begin{equation*}
\begin{split}
g(z)&\leq h(z)\sup_{\Omega\times\{|\tilde\xi|=1\}}\frac{g}{h}\leq h(z)(\sup_{\partial \Omega\times\{|\tilde\xi|=1\}}\frac{g}{h}+\frac{\mu}{\nu})\\
&\leq
h(z)(\beta^{\alpha}\sup_{\partial\Omega}u^{\vartheta}_{\gamma\gamma}+\hat{c}_{0}\beta^{-3}+\beta^{-8}).
\end{split}
\end{equation*}
By the virtue of (\ref{eq:defh}), (\ref{eq:defmu}) and
(\ref{eq:defnu}),
$$g(z)\leq
\frac{2\|\psi\|_{C^{1}}^{2}\beta^{\alpha}}{(\psi+\beta)^{\alpha}}\sup_{\partial\Omega}u^{\vartheta}_{\gamma\gamma}+
\frac{2\|\psi\|_{C^{1}}^{2}}{\beta^{\alpha}}(\hat{c}_{0}\beta^{-3}+\beta^{-8}),$$
where $\beta\in(0,\beta_{0})$, $\beta_{0}$ depends on
$\|\psi\|_{C^{2,1}},\,\|u^{\vartheta}\|_{C^{1}},\,\|f\|_{C^{1,1}}$, $C_{0}$
and $\|\rho\|_{C^{2,1}}$.
As $|\xi|^{2}+\xi_{0}^{2}=1$, we choose $\xi_{0}=0$,
$2\beta^{\alpha}=\varepsilon$ and
$C_{\varepsilon}=\frac{2c_{0}^{2}}{\beta^{\alpha}}(\hat{c}_{0}\beta^{-3}+\beta^{-8})$.
Thus, for any unit vector $\xi$,
\begin{equation}\label{eq:311'}
u^{\vartheta}_{\xi\xi}\leq\varepsilon\frac{\|\psi\|^{2}_{C^{1}(\Omega)}}{(\psi+\beta)^{\alpha}}\sup_{\partial\Omega}u^{\vartheta}_{\gamma\gamma}+C_{\varepsilon}.
\end{equation}
For general $\xi\in\mathbb{R}^{n}$, by timing $|\xi|^{2}$ with respect
to both sides of
(\ref{eq:311'}), we have
$$u^{\vartheta}_{\xi\xi}\leq
\varepsilon\frac{\|\psi\|^{2}_{C^{1}(\Omega)}}{(\psi+\beta)^{\alpha}}\sup_{\partial\Omega}u^{\vartheta}_{\gamma\gamma}+C_{\varepsilon}|\xi|^{2}.$$
That is (\ref{eq:weakin}), our proof is done.
\end{proof}

\bigskip
\section{Double Normal Derivative Estimate}\label{s5}

Now we prove Theorem \ref{th}. We need several steps. Our approach is still the
barrier.

Firstly, we study the second order derivative estimates on the boundary. For any point $x_{0}$ on the boundary, by a rotation and a translation, we can take $x_{0}$ as the origin. Then, choosing the principal
coordinates system at the origin, the boundary
$\partial\Omega$ is
represented by
$x_{n}=\rho(x')$ near the origin,
where $x'=(x_{1},\ldots,x_{n-1})$, $\rho$ is smooth as the smoothness
of $\partial\Omega$.

Let $T=(T_{i}^{j})$ be a skew-symmetric matrix, $\tau=(\tau_{1},\ldots,\tau_{n})$ be a vector field in $\Omega$ given by
$$\tau_{i}=T_{j}^{i}x_{j}+a_{i},\indent i=1,\ldots,n,$$
where $a_{i}$ is a constant.
We set
$$u^{\vartheta}_{(\tau)(\tau)}=:
  (u^{\vartheta}_{\tau})_{\tau}=:\tau_{i}\tau_{j}u^{\vartheta}_{ij}+(\tau_{i})_{j}\tau_{j}u^{\vartheta}_{i},$$
where $(u^{\vartheta}_{\tau})_{\tau}=\frac{\partial u_{\tau}}{\partial \tau}$, $(\tau_{i})_{j}=\frac{\partial\tau_{i}}{\partial x_{j}}$.

We have the following relationship.
\begin{proposition}[Lemma 2.1 in \cite{I}]
\begin{equation}\label{eq:dersp}
F[u^{\vartheta}_{\tau}]=(F[u^{\vartheta}])_{\tau},\indent
F^{ij}(u^{\vartheta}_{(\tau)(\tau)})_{ij}\geq(F[u^{\vartheta}])_{(\tau)(\tau)}.
\end{equation}
\end{proposition}

\begin{lemma}\label{lm:omega}
For any $\partial\Omega\ni x=(x_{1},\ldots,x_{n})=(x',x_{n})$ in the neighborhood of the origin, we have
\begin{equation}\label{eq:omega1}
u^{\vartheta}_{(\tau)(\tau)}(x)-u^{\vartheta}_{(\tau)(\tau)}(0)\leq \tilde{c}_{0}(|x'|^{2}+M|x'|^{4})\indent \text{in}\, B_{r_{0}}(0)\cap\partial\Omega,
\end{equation}
 where $\tau$ is the linear combination of $\eta^{i}$, $\tau=\alpha_{i}\eta^{i}$, $\sum \alpha^{2}_{i}=1$, $\eta^{i}$ is
 the annular vector
 field in the sense of the principal coordinates system at the origin,
\begin{equation*}
 \eta^{i}=\eta^{i}(x)=(1-\kappa_{i}(0)x_{n})\partial_{i}+\kappa_{i}(0)x_{i}\partial_{n},\indent i=1,\ldots,n-1,
\end{equation*}
and
 $$M=\sup_{x\in\partial\Omega}u^{\vartheta}_{\gamma\gamma}(x),$$
$\tilde{c}_{0}$ depends only on $\|\rho\|_{C^{2,1}}$, $\|u^{\vartheta}\|_{C^{1}(\overline\Omega)}$ and the bound of the second order mixed type derivatives, in particular, $\tilde{c}_{0}$ is independent with $\vartheta$.
\end{lemma}
\begin{proof}
Our proof will need several steps. Since
\begin{equation*}
\begin{split}
u^{\vartheta}_{(\tau)(\tau)}(x)&=\sum_{i,j,m,l}u^{\vartheta}_{ij}(\alpha_{m}\eta_{i}^{m})(\alpha_{l}\eta_{k}^{l})+u^{\vartheta}_{i}(\alpha_{m}\eta^{m}_{i})_{j}
(\alpha_{l}\eta_{j}^{l})\\
&=\left[\sum_{i,j,m}\alpha_{m}^{2}u^{\vartheta}_{ij}\eta^{m}_{i}\eta^{m}_{j}+\alpha^{2}_{m}u^{\vartheta}_{i}(\eta^{m}_{i})_{j}\eta^{m}_{j}\right]\\
&+\sum_{m\neq l}\sum_{i,j}[\alpha_{m}\alpha_{l}u^{\vartheta}_{ij}\eta^{m}_{i}\eta^{l}_{j}+\alpha_{m}\alpha_{l}u^{\vartheta}_{i}(\eta_{i}^{m})_{j}\eta^{l}_{j}]\\
&=P(x)+Q(x).
\end{split}
\end{equation*}
Setting
$$\omega=:u^{\vartheta}_{(\tau)(\tau)}(x)-u^{\vartheta}_{(\tau)(\tau)}(0).$$
Thus, we have
\begin{equation}\label{omec}
\omega(x)=P(x)-P(0)+Q(x)-Q(0).
\end{equation}
For some fixed $m$, let
\begin{equation*}
\omega^{m}(x)=:u^{\vartheta}_{(\eta^{m})(\eta^{m})}(x)-u^{\vartheta}_{(\eta^{m})(\eta^{m})}(0).
\end{equation*}
Thus, $$P(x)-P(0)=\sum_{m}\alpha^{2}_{m}\omega^{m}(x).$$

Firstly, we estimate the $P(x)-P(0)$. Recall the result in \cite{I}. For each $\omega^{m}(x)$, we have
$$\omega^{m}(x)\leq \tilde{c}_{0}(|x'|^{2}+M|x'|^{4}),$$
where $\tilde{c}_{0}$ depends only on $\|\rho\|_{C^{2,1}}$ and $\|u\|_{C^{1}}$.
Since $\sum\alpha_{m}^{2}=1$, we have
$$P(x)-P(0)\leq \tilde{c}_{0}(|x'|^{2}+M|x'|^{4}).$$

Then, by this result, we finish the proof of Lemma \ref{lm:omega}.
Recall \eqref{omec}. Thus, we need to estimate the left terms $Q(x)-Q(0)$. For fixed $m$ and $l$, where $m\neq l$, we take $\xi^{m}$ and $\xi^{l}$ as the projections of $\eta^{m}$ and $\eta^{l}$. We denote $\cos \theta^{l}=<\tilde\eta^{l},\tilde\xi^{l}>$, $\sin\theta^{l}=<\tilde\eta^{l},\gamma>$, where $\tilde\xi^{l}=\xi^{l}/|\xi^{l}|$, $\tilde\eta^{l}=\eta^{l}/|\eta^{l}|$.
Still, we have
\begin{equation*}
\begin{split}
&u^{\vartheta}_{(\eta^{m})(\eta^{l})}(x)-u^{\vartheta}_{(\eta^{m})(\eta^{l})}(0)\\
&\leq (1+4\|\rho\|_{C^{2,1}}|x'|^{2})(\varphi_{\tilde\xi^{m}\tilde\xi^{l}}(x)+(\varphi_{\gamma}(x)-u^{\vartheta}_{\gamma}(x))
\rho_{\tilde\xi^{m}\tilde\xi^{l}}(x))-u^{\vartheta}_{ml}(0)+H(x')\\
&\leq \varphi_{\tilde\xi^{m}\tilde\xi^{l}}(x)+(\varphi_{\gamma}(x)-u^{\vartheta}_{\gamma}(x))
\rho_{\tilde\xi^{m}\tilde\xi^{l}}(x)-\varphi_{ml}(0)+H(x')\\
&\leq H(x'),
\end{split}
\end{equation*}
where $H$ denotes a function satisfying
\begin{equation*}
H(x')\leq a|x'|+C|x'|^{2}+CM|x'|^{4}.
\end{equation*}
Thus,  by subtracting a linear function, we have
$$Q(x)-Q(0)\leq \tilde{c}_{0}(|x'|^{2}+M|x'|^{4}).$$
Accordingly, by \eqref{omec}, we have
$$\omega(x)\leq \tilde{c}_{0}(|x'|^{2}+M|x'|^{4}),$$
where $\tilde{c}_{0}$ may change from line to line, but all of them depend only on $\|\rho\|_{C^{2,1}}$, $\|u^{\vartheta}\|_{C^{1}(\overline\Omega)}$ and the bound of the second order mixed type derivatives, in particular $\tilde{c}_{0}$ is independent with $\vartheta$.
\end{proof}

Next we extend Lemma \ref{lm:omega} to the points in the
neighborhood of the origin. That
is
\begin{lemma}\label{lm:omega1}
For any $x\in\Omega$ near the origin
\begin{equation}\label{eq:gp0}
\omega(x)\leq C(|x'|^{2}+M|x'|^{4})+CM|x_{n}-\rho(x')|,
\end{equation}
where $C$ depends on $\|f\|_{C^{1,1}},r_{0}$, $C_{0}$ and
$\tilde{c}_{0}$, in particular, $C$ is independent with $\vartheta$.
\end{lemma}

\begin{proof}
Let
$\Omega_{r}=\{x\in\Omega|\rho(x')<x_{n}<\rho(x')+r^{4},|x'|<r^{4}\}$, and
$$v=C_{K}[(x_{n}-\rho(x'))^{2}+\beta(\rho(x')-x_{n})-\frac{\tilde{c}_{0}}{M}|x'|^{2}-\tilde{c}_{0}|x'|^{4}],$$
where $C_{K}\geq1$ is a constant. We choose $r$ as the highest
infinitesimal among
$r,\frac{1}{M},\beta$.
By the $(k-1)$-convexity of $\Omega$, $v$ is $k$-admissible. Thus, we
can choose $C_{K}$
which depends on $C_{0}$ so large that
$$F[v]\geq \delta_{0}+\sup_{\Omega}f^{1/k}$$
for some  positive constant $\delta_{0}\geq C_{0}$.
By the concavity of $F$ and Condition {\bf(H)}, one can verify
$$L[v]\geq F[u^{\vartheta}+v]-F[u^{\vartheta}]\geq F[u^{\vartheta}]-\sup_{\Omega}f^{1/k}\geq \delta_{0}\geq
C_{0}\geq -(F[D^{2}u^{\vartheta}])_{\tau\tau},$$
where $L[v]=F^{ij}\cdot v_{ij}$.

Then by (\ref{eq:dersp}), we have
\begin{equation}\label{eq:Mvome}
L(Mv)\geq-(F[D^{2}u^{\vartheta}])_{(\tau)(\tau)}\geq
-F^{ij}(u^{\vartheta}_{(\tau)(\tau)})_{ij}=L(-\omega).
\end{equation}
The boundary $\partial \Omega_{r}$ consists of three parts:
$\partial\Omega_{r}=\partial_{1}\Omega_{r}\cup\partial_{2}\Omega_{r}\cup\partial_{3}\Omega_{r}$,
where $\partial_{1}\Omega_{r}$ and $\partial_{2}\Omega_{r}$ are
respectively the graph
parts of $\rho$ and $\rho+r^{4}$,
$\partial_{3}\Omega_{r}$ is a portion of $\{|x'|=r\}$. $\partial_{1}\Omega_{r},\partial_{2}\Omega_{r}$ part, that is
$\rho(x')=x_{n}$ and
$x_{n}=\rho(x')+r^{4}$. We have
$$v\leq -\tilde{c}_{0}C_{K}(|x'|^{2}/M+|x'|^{4}).$$
$\partial_{3}\Omega_{r}$ part. $|x'|\leq r$, we have
$$v\leq
C_{K}[C_{1}r^{8}+C_{2}r^{4}-\frac{\tilde{c}_{0}}{M}r^{2}-\tilde{c}_{0}r^{4}]=-\frac{C_{K}\tilde{c}_{0}r^{2}}{M},$$
where $C_{1},C_{2}$ are constants depend on $\beta$ and $\Omega$.
Thus, we can choose $C_{K}$ which depends on $\|f\|_{C^{1,1}}$,
$C_{0}$ and $\Omega$ so
large that $v\leq-2$.

In sum, by (\ref{eq:omega1}), we have
$$Mv\leq-\omega\indent \text{on } \partial\Omega_{r}.$$
By the virtue of the classical weak maximum principle and
(\ref{eq:Mvome}), we have
$$Mv\leq -\omega \indent \text{in}\,\Omega_{r}.$$
That is (\ref{eq:gp0}).
\end{proof}

In order to find the barrier function, we still need a third order
derivative
estimate.
\begin{lemma}\label{lm:thirdd}
For any given $\sigma>0$ which is sufficiently small, we can find
a positive constant $C_{\sigma}$ such that
\begin{equation}\label{eq:unot}
(u^{\vartheta}_{(\tau)(\tau)})_{n}(0)\leq \sigma M+C_{\sigma},
\end{equation}
where $C_{\sigma}$ depends on $\tilde{c}_{0},\|f\|_{C^{1,1}}$, $C_{0}$
and the constant
$C_{\varepsilon,\delta}$ in \eqref{eq:weakin} with $\delta=\sigma^{4}$
and
$\varepsilon=\sigma^{8}$.
\end{lemma}
\begin{proof}
Let
\begin{equation*}
\tilde{v}(x)=\tilde{C}_{K}[(x_{n}-\rho(x'))^{2}+\beta(\rho(x')-x_{n})-\beta_{1}|x'|^{2}],
\end{equation*}
where $\beta,\beta_{1}>0$ are sufficiently small, $\tilde {C}_{K}>1$ is
sufficiently large.
Then for a sufficiently small $r>0$ which is different from which in
Lemma \ref{lm:omega1}.
$\tilde{v}$ is $k$-admissible in $\Omega_{r}$. We have
$$F[D^{2}\tilde{v}]\geq f_{0}\indent \text{in }\Omega_{r},$$
where $f_{0}\geq\delta_{0}+\sup_{\Omega}f^{1/k}$.

We want to explain that, if $r>0$ is sufficiently small and
$\tilde{C}_{K}>1$ is
sufficiently large,
$(rM+C_{r})\tilde{v}$ is a sub-barrier of $\omega$, where
$C_{r}=r^{-4}C_{\varepsilon,\delta}$. We have
\begin{equation}\label{eq:unobound}
(rM+C_{r})\tilde{v}\leq -\omega\indent\text{on }\partial\Omega_{r}.
\end{equation}
Indeed, on $\partial_{1}\Omega_{r}$, we have
$$\tilde{v}(x)\leq -\beta_{1}\tilde{C}_{K}|x'|^{2},$$
and hence
\begin{equation*}
\begin{split}
(rM+C_{r})\tilde{v}&\leq
(rM+r^{-4}C_{\varepsilon,\delta})(-\beta_{1}\tilde{C}_{K}|x'|^{2})\\
&\leq
-\beta_{1}\tilde{C}_{K}Mr|x'|^{2}-\beta_{1}\tilde{C}_{K}C_{\varepsilon,\delta}\frac{|x'|^{2}}{r^{4}}\\
&\leq-\omega
\end{split}
\end{equation*}
provided $\tilde{C}_{K}$ is
sufficiently large.
\newline On $\partial_{2}\Omega_{r}$, $\tilde{v}<-\frac{1}{2}\beta
\tilde{C}_{K}r^{4}$. By
Lemma \ref{lm:weakin}
$$-\omega(x)\geq-r^{8}M-C_{r^{8},r^{4}}\geq
-r^{8}-\frac{C_{\varepsilon,\delta}}{r^{8}}.$$
Since $r=o(\beta)$, (\ref{eq:unobound}) holds.
\newline On $\partial_{3}\Omega_{r}$, we have
$$\tilde{v}=\tilde{C}_{K}[(x_{n}-\rho(x')-\frac{\beta}{2})^{2}-\frac{\beta^{2}}{4}-\beta_{1}|x'|^{2}]<
-\beta_{1}\tilde{C}_{K}r^{2}.$$
By (\ref{eq:gp0}),
$$-\omega\geq -C(1+r^{4}M),$$
hence (\ref{eq:unobound}) holds. In $\Omega_{r}$, we have
$L(\tilde{v})\geq \delta_{0}$.
We still use the approach as before.
One can verify that
$$L[(rM+C_{r})\tilde{v}]\geq C_{0}.$$
 Thus, by the virtue of  (\ref{eq:dersp}) and Condition {\bf(H)}
$$L[(rM+C_{r})\tilde{v}]\geq-L(\omega)\indent \text{in }\Omega_{r}.$$
In sum,
$$(rM+C_{r})\tilde{v}\leq-\omega\indent \text{in }\Omega_{r}.$$
Thus, we have
$$(u^{\vartheta}_{(\tau)(\tau)})_{n}(0)\leq (rM+C_{r})|\tilde{v}_{n}(0)|\leq
C_{\tilde{v}_{n}}(rM+C_{r}).$$
Then, let $\sigma=C_{\tilde{v}_{n}}r$, that is (\ref{eq:unot}).
\end{proof}
Finally, we finish the proof of Theorem \ref{th}. Thus, we complete the
whole proof.
\begin{proof}[{\bf{End of the Proof to Theorem \ref{th}}}]
For any boundary point $x_{0}$ and arbitrary tangential unit vector field $\xi$ on
$\partial\Omega$
near $x_{0}$. We set $x_{0}$ to be the origin, Then, $\partial\Omega$
is locally given by
$x_{n}=\rho(x')$, $\tau$ and $\eta^{m}$ are defined as which in Lemma \ref{lm:omega}. At the origin,
$$\tau(0)=\alpha_{m}\partial_{m}.$$
Since $\sum_{m}\alpha^{2}_{m}=1$,  $\tau(0)$ is a tangent vector. Thus, we can choose some appropriate $\alpha_{m}$ such that $\xi=\tau$ at $x_{0}$.
One can verify that at $x_{0}$
$$u^{\vartheta}_{n(\xi)(\xi)}=u^{\vartheta}_{n(\tau)(\tau)}=u^{\vartheta}_{nij}\tau_{i}\tau_{j}+(\tau_{i})_{j}\tau_{j}u^{\vartheta}_{ni},$$
$$(u^{\vartheta}_{(\tau)(\tau)})_{n}=u^{\vartheta}_{ijn}\tau_{i}\tau_{j}
+u^{\vartheta}_{ij}(\tau_{i}\tau_{j})_{n}+((\tau_{i})_{j}\tau_{j})_{n}u^{\vartheta}_{i}+(\tau_{i})_{j}\tau_{j}u^{\vartheta}_{in}.$$
At the origin, $\tau_{n}=0$, $u_{ij}(0)$ is bounded when $i,j=1,\ldots,n-1$. Then
$$|u^{\vartheta}_{ij}(\tau_{i}\tau_{j})_{n}+((\tau_{i})_{j}\tau_{j})_{n}u^{\vartheta}_{i}|\leq C_{\sigma}.$$
Thus, we have
$$u^{\vartheta}_{n(\xi)(\xi)}(0)=(u^{\vartheta}_{(\tau)(\tau)})_{n}-(u^{\vartheta}_{ij}(\tau_{i}\tau_{j})_{n}+((\tau_{i})_{j}\tau_{j})_{n}u^{\vartheta}_{i})
\leq \sigma M+C_{\sigma}.$$
Since this result is independent with the choice of $x_{0}$, we have, for any tangent vector field $\xi=\sum_{i=1}^{n-1}\alpha_{i}(\partial_{i}+\rho_{i}(x')\partial_{n})$, $\sum^{n-1}_{i=1}\alpha^{2}_{i}=1$
\begin{equation}\label{eq:taylor2}
u^{\vartheta}_{n(\xi)(\xi)}(x)\leq \sigma M+C_{\sigma}\indent \text{on}\,\partial\Omega.
\end{equation}

Choosing a new coordinates system, we suppose the maximum $M$ is
attained at the origin $0\in
\partial\Omega$. Then near the origin, we set $G(x')=u^{\vartheta}_{n}(x',\rho(x'))$ defined on $\partial\Omega$, by the Taylor expansion on the boundary, for $h=(h_{1},\ldots,h_{n})=(h',h_{n})\in \partial\Omega$ near the origin
\begin{equation*}
G(h')=G(0)+\sum_{i=1}^{n-1}h_{i}\partial_{i}G(0)+\int_{0}^{1}(1-t)(\sum_{i=1}^{n-1}h_{i}\partial_{i})^{2}G(th')dt.
\end{equation*}
One can verify that $$G(0)=u^{\vartheta}_{n}(0),$$
$$\sum_{i=1}^{n-1}h_{i}\partial_{i}G(x')=\sum_{i=1}^{n-1}h_{i}(\partial_{i}+\rho_{i}(x')\partial_{n})u^{\vartheta}_{n}(x',\rho(x')).$$
If we choose $\alpha_{i}=h_{i}/|h'|$, where $|h'|=(\sum_{m=1}^{n-1}h^{2}_{m})^{1/2}$, we have
$$\sum_{i=1}^{n-1}h_{i}\partial_{i}G(x')=|h'|u^{\vartheta}_{n\xi}(x',\rho(x')).$$
Still, we have
$$(\sum_{i=1}^{n-1}h_{i}\partial_{i})^{2}G(th')=u^{\vartheta}_{n(\xi)(\xi)}(th',\rho(th'))|h'|^{2}.$$
Then, by \eqref{eq:taylor2}, we have
$$u^{\vartheta}_{n}(h)\leq u^{\vartheta}_{n}(0)+a|h'|+(\sigma M+C_{\sigma})|h'|^{2},$$
where $a$ is the bound of the mixed second order derivatives. By
subtracting a linear
function and the gradient estimate, we can obtain
$$L[(\sigma M+C_{\sigma})\tilde{v}]\geq L(-u^{\vartheta}_{n})\indent \text{in
}\Omega\cap B_{r},$$
$$(\sigma
M+C_{\sigma})\tilde{v}\leq-(u^{\vartheta}_{n}(x)-u^{\vartheta}_{n}(0))\indent\text{on
}\partial(\Omega\cap B_{r}).$$
Thus, we have
$$-(u^{\vartheta}_{n}(x)-u^{\vartheta}_{n}(0))\geq (\sigma M+C_{\sigma})\tilde{v}(x)\indent
\text{in }\Omega\cap
B_{r},$$
namely,
$$M=u^{\vartheta}_{nn}(0)\leq -\tilde{v}_{n}(0)(\sigma M +C_{\sigma})\leq
\tilde{C}_{K}\sigma(\sigma
M+C_{\sigma}).$$
Let $\sigma^{2}<1/2$, one can verify that
$$M\leq 2\tilde{C}_{K}\sigma C_{\sigma}\leq
\frac{C_{\sigma}}{\sigma}$$
holds for any $\sigma\in(0,\sigma_{0})$. From the computation above
$$\frac{C_{\sigma}}{\sigma}=O(\sigma^{-m}),$$
where $m$ is a positive constant. Thus,
$$\sup_{\partial\Omega}u^{\vartheta}_{\gamma\gamma}\leq
\frac{C_{\sigma_{0}}}{\sigma_{0}},$$
where $\sigma_{0}$ depends only on $\|\psi\|_{C^{3}(\Omega)},\Omega$,
$C_{0}$,
$\|\rho\|_{C^{2,1}}$ and $\|\varphi\|_{C^{3,1}}$,
in particular, $C_{\sigma_{0}}$ and $\sigma_{0}$ are independent with
$\vartheta$.
Proof is done.
\end{proof}

\vspace{5mm}

{\section*{\small{Acknowledgements}}
{\small{This work was supported by the Fundamental Research Funds for
the Central Universities No. 2012201020206 and
the National Science Foundation of China No. 11171261. The authors would like to appreciate
Guji Tian for his helpful conversations.}}

\end{document}